\documentclass[11pt]{article}

\usepackage[a4paper, margin=2.5cm]{geometry}
\usepackage{amsmath, amssymb, amsthm,dsfont}
\usepackage{graphicx}

\usepackage{hyperref}
\hypersetup{
    colorlinks,
    linkcolor={red!50!black},
    citecolor={blue!50!black},
    urlcolor={blue!80!black}
}
\newtheorem{lemma}{Lemma}

\numberwithin{equation}{section}
\usepackage{algorithm}
\usepackage{algpseudocode}
\usepackage{tikz}
\newtheorem{remark}{Remark}
\usepackage{comment}
\usepackage{xcolor}
\usepackage{todonotes}

\usepackage{hyperref}
 [
     colorlinks=true, 
     linkcolor=blue!60!black, 
     citecolor=green!60!black, 
     filecolor=magenta, 
     urlcolor=cyan!60!black
 ]
\usepackage{pgfplots}
\pgfplotsset{compat=1.18} % ou la version installée
\usepackage{graphicx}
\pgfplotsset{
  table/search path={PINNs-example1,PINNs-example2,DR-example2,DR-example1,PINNs-example3}
}
\usepackage{float}
\title{Goal oriented error estimation for adaptive sampling of PINNS }

\author{
    Medard Govoeyi$^{1}$, Thomas Richter$^{2}$, \\
    Otto-von Guericke University,
    \texttt{medard.govoeyi@ovgu.de, thomas.rchter@ovgu.de } 
}
\date{}

\begin{document}

\maketitle
\begin{abstract}
Physics-Informed Neural Networks (PINNs) are mesh-free approaches for the numerical approximation of partial differential equations, where a neural network is trained by minimizing a loss function derived from the governing equations and boundary conditions. The Deep Ritz method can be interpreted as a particular variational form of a PINN, where the loss corresponds to the minimization of an energy functional associated with a symmetric positive definite problem.

In this work, we study the approximation of the Laplace equation using both the classical PINN formulation and its variational counterpart, the Deep Ritz method, with the objective of accurately estimating prescribed goal functionals. When standard sampling strategies, such as uniform or loss-based sampling, are employed during training, the convergence of the functional error and the attained minimal functional value can be slow.

To address this issue, we introduce a functional-oriented importance sampling strategy that can be applied to both PINNs and the Deep Ritz method. The key ingredient is the construction of a reliable and accurate estimator for the error in a given quantity of interest. This estimator is derived using concepts from the Dual Weighted Residual (DWR) framework and is implemented entirely within the neural network setting. It is then used to adaptively guide the sampling of training points in the computational domain, focusing computational effort on regions that have the strongest influence on the functional value.

Numerical experiments demonstrate that the proposed adaptive sampling strategy significantly accelerates the convergence of the functional error and improves the minimization of the target functional during training for both PINN and Deep Ritz formulations.
\end{abstract}

% \includegraphics{figs/cas-grabs.pdf}

% \begin{highlights}
% \item Research highlights item 1
% \item Research highlights item 2
% \item Research highlights item 3
% \end{highlights}

\textbf{keywords}
Scientific Machine Learning ,Physics-Informed Neural Networks (PINNs),
a posteriori error estimation,
adaptive sampling

\section{Introduction}

The emergence  of machine learning in the last decade has contributed to develop multiple tools to solve differential equations. One class of such methods is the area of  Physics Informed Neural Networks (PINNs), where a neural network is trained with a loss function that penalizes the residual of the differential equation, such that no (or little) data is needed. It's Lagaris et al. who presented this approach for the first time  in \cite{Lag98}. Often, PINNs formulate loss functions that are composed of different parts to match the differential equation, the boundary conditions, the residual of the differential equation  or further constraints. Many contributions and developments have been recently published in the field of PINNs, notably in \cite{Rai17}, \cite{Rai19}, where PINNs are used to approximate complex differential equations and Inverse problems. PINNs considers the solution of the differential equations as neural network and train it to satisfy the differential equations and the boundary conditions. Unlike traditional methods for solving differential equations, Physics-Informed Neural Networks (PINNs) do not require the construction of a mesh. As a result, the geometry of the domain is not a limiting factor. This key advantage enables PINNs to handle high-dimensional partial differential equations (PDEs). In \cite{Hu24}, this method was successfully applied to solve PDEs in dimension 
$d=100\,000$ and in \cite{Hu24-1,Zha24}, PINNs have been employed in high-dimensional settings, with the aim of overcoming the curse of dimensionality.

One simple method from the PINNs  family is the Deep Ritz method, introduced by Li and E~\cite{E18}. It directly minimizes the energy functional of a symmetric positive definite variational form, expressing the solution by a neural network. Deep Ritz is attractive as it is open to analytical approaches that close to established numerical analysis~\cite{MuellerZeinhofer2021,MR}. 

In our contribution we will also focus on Deep Ritz, as its very close vicinity to energy minimization methods opens a path to mathematical analysis and gives a direct connection to established numerical approaches. The aim of this paper is to explore these connections and to leverage classical concepts such as error estimation and adaptivity to improve the efficiency.

PINNs show potential as an alternative to classical discretization methods, however they suffer  in terms of lacking robustness in their approximation quality and in particular they have to undergo the training process. For complex problems, training the neural network, which means solving a highly unstructured optimization problem, is still the challenging part. Compared to the finite element method, where a solution is obtained by solving a linear (or nonlinear) system of equations, which can be done in optimal complexity, PINNs require the approximation of a large and non convex optimization problem.

\paragraph{Adaptive Sampling Approaches}
To optimize the training process, different adaptive methods have been developed to improve training.
%+These methods use the loss to sample the training data (collocations points)   \cite{Sub22,Tor24,Roy24,Liu23}.
We can distinguish two classes of methods: Residual-based Adaptive Refinement (RAR) and Residual-based Adaptive Distribution (RAD).
The objective of RAR is to add, during training, new collocation points that have the largest residuals, selected from a pool of candidate points \cite{Lu2021DeepXDE}.
RAD, on the other hand, considers a probability density function (PDF) derived from the residual, and periodically replaces the current collocation points with new points sampled using this PDF \cite{Nabian2021, wu2023comprehensive, Gao2023FIPINNs, Tang2023DAS}.
There is also RAR-D, which combines both approaches. It consists of adding new points as in RAR, but these new points are obtained by sampling using the residual-based PDF distribution \cite{wu2023comprehensive}.
All these methods are based on the pointwise evaluation of the residual that enters the loss. 

In~\cite{Lu2026} the authors provide an error estimator that is then used to sample the collocation points. Their estimator builds on classical gradient recovery strategies of Zienkiewicz \& Zhu type~\cite{ZZ1,ZZ2}. In particular for solutions with localized features they demonstrate excellent performance improvement by using this error-informed adaptive procedure. One drawback of the recovery method is the use of an underlying mesh structure, which will not be feasible for very high dimensional problems, whose efficient approximation is one of the theoretical advantages of PINNs.

\paragraph{Aims and scope} 
In this paper, we develop an adaptive training method optimized to yield functional outputs of the solution. 
We aim at steering the training process such that the solution measured in $J(u)$, where $J$ is the error functional, is trained as fast and as accurate as possible. Like the authors of~\cite{Lu2026} we draw from ideas from classical finite element analysis and base the estimator on the Dual Weighted Residual method (DWR) which has been introduced by Becker and Rannacher~\cite{BeckerRannacher1995,BeckerRannacher2001} and developed since then into a standard tool in finite element analysis~\cite{Endtmayer2024,Becker2025}. 
DWR is a computable error approximation (not a strict bound) that measures the error in arbitrary output functions, which can be norms, point values, averages, etc. DWR introduces a dual equation for obtaining sensitivities of the solution with respect to the functional. The DWR estimator delivers bounds on the approximation error of the PINN but furthermore, the DWR error can be localized and used to adaptively control the discretization, see~\cite{RichterWick2015}. We will exploit this localization for adaptive sampling of collocation points used in training the PINN.  

The DWR estimator has been extended to estimate the Deep Ritz error as highly accurate stopping criteria during training~\cite{MR}. Similar to~\cite{Lu2026}, the limitation of this approach is the use of a mesh-structure needed for evaluating the estimator.
However, the key finding from~\cite{MR} is the combination of machine learning techniques with finite element analysis techniques such as DWR method. Furthermore, PINN approximation have also been used to obtain sensitivity information within a finite element application of DWR~\cite{RothSchroederWick2022}. This attempt can be considered as the opposite approach to~\cite{MR}, where the solution itself is a PINN but the adjoint is approximated using finite elements. These two paper numerically demonstrated that one typical problem of DWR based estimation appears to be less troublesome: While it is usually important to approximate the dual weights in function spaces that are substantially richer than the discrete solution space, which often results in the effort required to estimate the error being greater than the effort for solving the problem, the missing linear structure of neural network sets prevents this problem from occurring. 

Here, we will introduce a fully PINN-based approach for error estimation that does not require any underlying mesh structure for the evaluation of the estimator. Furthermore, we will present a localization of the estimator to be used for sampling collocation points. 

\paragraph{Overview}
In Section \ref{section2}, we present the Deep Ritz Method.
Section \ref{section3} is devoted to the Dual Weighted Residual (DWR) method, which we first introduce in the context of the finite element method, before transferring it to PINNs. Section~\ref{section4} describes our approach for the adaptive control of PINNs by choosing collocation points guided by the error estimator.
Finally, in Section \ref{section5}, we explain how this estimator is used to enhance the learning process and we report several numerical test cases.

\section{Physics Informed Neural Networks}  \label{section2}
  
We briefly introduce the basic concepts of PINNs \cite{Rai19} as neural network based approaches for the approximation of PDEs. The general idea is that the solution is represented by a neural network function $u_\theta:\mathbb{R}^d\to\mathbb{R}^c$ with trainable parameters $\theta$ that are obtained by a minimization problem with a loss function which is determined by the PDE itself, therefore the label ``physics informed''.

Considering the Laplace problem $-\Delta u = f$ in $\Omega$ with $u=g$ on the boundary $\partial\Omega$, the most simple type of PINN is a penalization of the strong residual accompanied by a penalization of the boundary condition
\begin{equation}\label{loss:PINN}
E_\lambda(u_\theta) := \|-\Delta u_\theta - f\|^2_\Omega + \lambda \|u_\theta-g\|^2_{\partial\Omega}
\end{equation}
with a tunable parameter $\lambda>0$. 
The corresponding loss function is given by \cite{Rai19}
\begin{equation}\label{PINN:loss}
E_{\lambda, M, N} := \frac{1}{N}\sum_{i=1}^{N} \left| -\Delta u_\theta (x_i) - f(x_i) \right| ^2
+  \frac{\lambda}{M} \sum_{j=1}^{M}\left|u_\theta(y_j)-g(y_j)\right|^2, 
\end{equation}
where $x_i \in \Omega$ and $y_j \in \partial \Omega$ for $i=1,\dots,N$ and $j=1,\dots,M$ are the collocation points within the domain and on the boundary.
\paragraph{Deep Ritz:} 
The Deep Ritz method \cite{E18} is a variant for symmetric positive definite problems that are given as minimizer of an energy form
\begin{equation*}
    u_\theta \in V, \quad  E(u) \leq E(v), \quad \forall v \in V,
\end{equation*}     
which, for the Poisson problem $-\Delta u=f$ just gets
\[
E(v) = \frac{1}{2}\int_\Omega |\nabla v(x)|^2\,\text{d}x - \int_\Omega f(x)v(x)\,\text{d}x,
\]
As it is difficult to strongly integrate boundary values into neural networks, the loss function again is enriched by a boundary penalization term, such that Deep Ritz is based on minimizing the loss
\begin{equation}\label{loss:deepritz}
E_\lambda(u_\theta) = \frac{1}{2}\int_\Omega |\nabla u_\theta(x)|^2\,\text{d}x - \int_\Omega f(x)u_\theta(x)\,\text{d}x + \frac{\lambda}{2} \|u_\theta-g\|^2_{\partial \Omega}.
\end{equation}
An obvious benefit of Deep Ritz over the general PINN approach is that only first derivatives are used. Minimizing the Energy function with boundary penalization in $H^1(\Omega)$ yields the disturbed solution
\begin{equation}\label{laplace1:dist}
\begin{aligned}
    - \Delta u_\lambda & = f\text{ in } \Omega \\
        u_\lambda + \frac{1}{\lambda}\partial_n u_\lambda  & =0\text{ on }\partial \Omega 
\end{aligned}
\end{equation}
Hereby we introduce a model error $(u-u_\lambda)$ which is of order $O(\lambda^{-1})$ and which has been analyzed in~\cite{Duan} and~\cite{MR}.

For the application to neural network functions $u_\theta\in V_{\cal N}$ we must ensure that the set of network functions $V_{\cal N}$ is conforming with the losses~\eqref{loss:PINN} and \eqref{loss:deepritz} which means $V_{\cal N}\subset H^2(\Omega)$ for the general PINN and $V_{\cal N}\subset H^1(\Omega)$ for Deep Ritz. However, as there are no efficient computational ways to analytically evaluate the integrals and norms in~\eqref{loss:PINN} and~\eqref{loss:deepritz}, numerical quadrature must be used such that $V_{\cal N}\subset C(\bar\Omega)$ will always be required. In general we will use neural network architectures that guarantee $V_{\cal N}\subset C^\infty(\bar\Omega)$.
%%
%\begin{equation*}
%  V_\mathcal{N} = \left\{ u_\theta\in C^\infty(\bar\Omega),\;  \theta \in \Theta  \right\}, 
%\end{equation*}
%
%where $\Theta \subset\mathbb{R}^{n_P}$ is the space of admissible parameters. 
Details on the design of neural networks follow in Section~\ref{subsection2.2}. 

The Deep Ritz method uses a Monte Carlo  approximation of loss for training
\begin{equation}\label{DR:loss}
    E_{\lambda, M, N}(u_\theta) := \frac{\vert \Omega \rvert }{N}\sum_{i=1}^{N} \frac{1}{2} \lvert \nabla u_\theta (x_i) \rvert ^2 - f(x_i) u_\theta(x_i)
    + \lambda \frac{\vert \partial \Omega \rvert }{2M} \sum_{j=1}^{M}|u_\theta(y_j)-g(y_j)|^2 ,   
\end{equation}
with $x_i \in \Omega$ and $y_j \in \partial \Omega$ for $i=1,\dots,N$ and $j=1,\dots,M$.

The Deep Ritz approximation of the solution of the Laplace  equation is the the solution obtained by training the neural network $u_\theta\in V_\mathcal{N}$ using $E_{\lambda, N, M}(\cdot)$ as the loss function.  We denoted the solution $u_{\lambda, \mathcal{N}}$.

\subsection{Neural Network Architectures}  
As neural network architecture we use standard Feed Forward Neural Networks, also denoted as multilayer perceptrons (MLP), recursively defined as
 \begin{equation}\label{mlp}
     \begin{aligned}
         u_\mathcal{N}^{(0)} (x) &:=x ,\\
         u_\mathcal{N}^{(l)} (x) &: = \sigma \left( W^lu_{\mathcal{N}}^{(l-1)}(x)+ b^l \right)\text{ for } l=1,2, \dots,L,\\
           u_\mathcal{N}(x) & = W^{L+1} v_\mathcal{N}^{(L)} + b^{L+1},
     \end{aligned}
 \end{equation}
 where weight matrices $W^l \in \mathbb{R}^{N_l \times N_{l-1}}$ and bias $b^l \in \mathbb{R}^{N_l}$ are the trainable parameters.
Each $u_\mathcal{N}^l$ represents the $l$-th layer, $N_l$ is the number of neurons of the $l$-th layer and $L$ is the total number of the layers. By $\sigma$ we denote the activation function which is acting component-wise on the vectors $(W^lu_{\cal N}^{(l-1)}+b^l)$. Given $\sigma\in C^\infty(\mathbb{R})$, the network, as composition of smooth activation functions and affine functions, has the same regularity. The often used activation functions $\tanh$ has this regularity, while the ReLU activation is continuous but only piecewise smooth. A smooth variant of ReLU is the GeLU activation~\cite{GELU2023} which is also well suited for PINNs.
 
In literature, superior performance of PINNs (including Deep Ritz) is often reported for so-called residual neural networks~\cite{E18}  that impose skip-connections, meaning that the inner layers can be replaced by
  \begin{equation}\label{resnet}
    u_\mathcal{N}^{(l)} (x) : =u_{\cal N}^{(l-1)}(x)+ \sigma \left( W^lu_{\mathcal{N}}^{(l-1)}(x)+ b^l \right).
  \end{equation}
  The set of neural network functions is given by
  \begin{equation}
     V_\mathcal{N} = \big\{ u_\mathcal{N}\text{ given by }\eqref{mlp}\text{ or }\eqref{resnet}\text{ with } W^l \in \mathbb{R}^{N_l \times N_{l-1}}, b^l \in \mathbb{R}^{N_l},\; l=1,\dots,L+1 \big\}. 
  \end{equation}
  Depending on the number and width of the layers, the set $V_{\cal N}$ is high dimensional with 
  \[
  \#V_{\cal N} = \sum_{l=1}^{L+1}N_l(1+N_{l-1})
  \]
  free parameters. 
  
\subsection[Deep Ritz as variational problem and the role of the neural network set V\_N]{Deep Ritz as variational problem and the role of the neural network set $V_{\cal N}$}\label{subsection2.2}

For the following we define the $L^2$-inner products on $\Omega$ and the boundary $\partial\Omega$
 \begin{equation*}
         \langle u, v \rangle_\Omega := \int_\Omega u(x) v(x) dx,\quad
         \langle u, v \rangle_{\partial \Omega} := \int_{\partial \Omega} u(x) v(x) dx.
 \end{equation*}

If we minimize the Deep Ritz energy functional~\eqref{loss:deepritz} over a Hilbert space $V$, e.g. $V=H^1(\Omega)$, the minimum $u_\lambda\in V$ is uniquely determined as solution to the variational problem
\[
E_\lambda(u_\lambda)\le E(v)\;  \forall v \in V \quad\Leftrightarrow \quad
\langle \nabla u_\lambda,\nabla\phi\rangle_\Omega 
+\lambda \langle u_\lambda-g,\phi\rangle_{\partial\Omega}
=\langle f,\phi\rangle_\Omega, \quad  \forall  v \in V.  
\] Necessary for this equivalence is the linear structure of the space $V$, which, in general, is not given for the set of neural network functions. For example, the set
\[
V_{\cal N}=\{ x \mapsto w_2\tanh(w_1x+b_1)+b_2:\; w_1,w_2,b_1,b_2\in\mathbb{R}\}
\]
is not closed with respect to the addition as, for instance it holds
\[
\tanh(x)+\tanh(x+1)\not\in V_{\cal N}
\]
although $\tanh(x),\tanh(x+1)\in V_{\cal N}$. 
Nevertheless, $V_{\cal N}\subset H^1(\Omega)$, and, therefore, it holds
\[
\min_{v\in V_{\cal N}}E_\lambda(v)\ge \min_{v\in H^1(\Omega)}E(v)> -\infty
\] 
and since $V_{\cal N}$ is finite dimensional there exists a (not necessarily unique) minimizer $u_\theta\in V_{\cal N}$ with
\begin{equation}\label{Eq-min-energy}
E_\lambda(u_\theta) = \min_{v\in V_{\cal N}}E_\lambda(v).
\end{equation}

\begin{lemma}[Variational Structure of Deep Ritz]\label{lemma:var}
Let $u_\theta \in V_{\mathcal N}$ be parameterized by
\[
\theta \in \Theta = (W_1,\dots,W_{L+1},b_1,\dots,b_{L+1}) = \mathbb{R}^{n_p}.
\]
Since the parameter space $\Theta$ has a linear structure, any (not necessarily unique) minimizer $u_\theta $ of \eqref{Eq-min-energy} can be characterized as follows
$$
\frac{d}{ds}E(u_{\theta+s\beta})\Big|_{s=0}\stackrel{!}{=}0\quad\forall \beta\in\Theta.
$$

In terms of neural network functions, this characterization corresponds to the Euler-Lagrange equation 
\begin{equation}\label{Euler-Lagrange}
u_\theta\in V_{\cal N}:\quad      
\langle \nabla u_\theta, \nabla v_\beta' \rangle_\Omega + \lambda \langle u_\theta-g , v_\beta' \rangle_{\partial \Omega} = \langle f,  v_\beta' \rangle_\Omega\quad  \forall v_\beta' \in   V_\mathcal{N}^\prime(u_\theta),
\end{equation}
where 
\begin{equation*}
    V_\mathcal{N}^\prime(u_\theta) := \left\{ x \mapsto \frac{d}{ds} u_{\theta+s\beta}(x)\big|_{s=0}:\quad \beta \in \Theta=\mathbb{R}^{n_p}\right\}.
\end{equation*}
In general $V_{\cal N}\neq V_{\cal N}'(u_\theta)$.
\end{lemma}

 \paragraph{Example 1}\label{Example}
 Let us consider the shallow neural network,
 \begin{equation*}
 V_{\cal N}:=\{u_\theta: x \mapsto \theta_1 \sigma( \theta_2 x + \theta_3 ) + \theta_4,\quad \theta_1,\dots,\theta_4\in \mathbb{R}\}
 \end{equation*}
 with a differentiable activation function $\sigma(\cdot)$. Then, the set of derivative networks is given by
 \[
 V_{\cal N}'(u_\theta):=
 \{
 x \mapsto \beta_1 \sigma( \theta_2 x + \theta_3 ) + \beta_4
 +
 \theta_1 \sigma'( \theta_2 x + \theta_3 )
 (\beta_2 x+\beta_3),\quad
 \beta_1,\dots,\beta_4\in \mathbb{R}
 \}
 \]
 A subset of this derivative set is given, if only the most outer weights $\beta_1$ and $\beta_4$ are considered and $\beta_2$ and $\beta_3$ are chosen as zero
  \[
 \tilde V_{\cal N}'(u_\theta):=
 \{
 \beta_1 \sigma( \theta_2 x + \theta_3 ) + \beta_4,\quad
 \beta_1,\beta_4\in \mathbb{R}
 \}\subset V_{\cal N}'(\theta).
 \]
 This subset of network functions is easily realized by keeping the architecture of $V_{\cal N}$, freezing the inner coefficients $\theta_2,\theta_3$ and only leaving the outer weights as tunable parameters.

The role of the derivative set $V'_{\cal N}(u_\theta)$ and its approximation $\tilde V'_{\cal N}$ lies in the interpretation of the Euler-Lagrange equation~\eqref{Euler-Lagrange} as Petrov-Galerkin formulation with different trial and test sets. This interpretation however is only rather formal, as the test-set $V_{\cal N}'(u_\theta)$ depends on the solution $u_\theta$ and its determining weights $\theta\in \Theta$. We will exploit this structure for defining proper adjoint spaces when it comes to error estimation. 

\section{A posteriori error estimates for the Deep Ritz method}\label{section3}  

Based on~\cite{MR}, we introduce an a posteriori error estimator for functional values  the error $u-u_{\lambda, \mathcal{N}}$. In contrast to~\cite{MR} we will establish a practical error estimator that is completely based on neural network approximations and does not utilize finite element solutions for approximating the adjoint solution. Furthermore, we will derive a localization of the error estimator for adaptively steering the training process of the neural network solution.

\subsection{The Dual Weighted Residual method for the finite element method}

We start with a brief introduction of the classical DWR method for elliptic problems, see~\cite{BeckerRannacher1995,BeckerRannacher2001}.
Let $a:H^1(\Omega)\times H^1(\Omega)\to\mathbb{R}$ be the bilinear form of the weak formulation of the partial differential equation, i.e.
\[
u\in H^1(\Omega):\quad 
a(u,\phi)=F(\phi)\quad\forall \phi\in H^1(\Omega),
\]
which, for the Poisson problem with boundary penalization gets
\[
a(u,\phi)=\langle\nabla u,\nabla\phi\rangle_{\Omega}
+
\lambda\langle\nabla u,\nabla\phi\rangle_{\Omega},\quad
F(\phi) = \langle f,\phi\rangle_{\Omega}.
\]
Next, let $J:H^1(\Omega)\to\mathbb{R}$ be a linear error functional and for  our purpose it is sufficient to assume that $J$ can be written in the form
\[
J(\phi)=\langle j,\phi\rangle_{\Omega}
\]
with $j\in L^2(\Omega)$. Hence, $J(\phi)$ can be considered a weighted average of $\phi$. We introduce the adjoint equation for $z\in H^1(\Omega)$
\begin{equation}\label{adjoint:var}
a(\phi,z)=J(\phi)\quad\forall \phi\in H^1(\Omega).
\end{equation}
Then, if we have any approximation $\tilde u\in \tilde V\subset H^1(\Omega)$ the error representation
\begin{equation}\label{errorrepresentation}
	J(u-\tilde u) = a(u-\tilde u,z) = \langle f,z\rangle_{\Omega}-a(\tilde u,z)
\end{equation}
holds. Instead of the true solution $u\in H^1(\Omega)$ the adjoint solution $z\in H^1(\Omega)$ is required for its evaluation. 

\begin{remark}[Approximation of the adjoint solution]
	In the context of the finite element method $u_h\in V_h\subset H^1(\Omega)$ comes from a discrete subspace. Then, it is not sufficient to approximate the adjoint in the same $z_h\in V_h$ by
	\[
	a(\phi_h,z_h) = J(\phi_h)\quad\forall \phi_h\in V_h
	\]
	since Galerkin orthogonality would counteract any representative content of the error identity~\eqref{errorrepresentation}:
	\[
	J(u-u_h) = \langle f,z\rangle_{\Omega}-a(u_h,z)\neq
	\langle f,z_h\rangle_{\Omega}-a(u_h,z_h)=0.
	\]
	Instead, one usually refers to higher order approximations or higher order reconstructions $u_{h}^{(*)}\in V_h^{(*)}$ in enriched spaces $V_h^{(*)}\subset H^1(\Omega)$ with $V_h^{(*)}\not\subset V_h$. Global solutions in finer spaces should be avoided since this would mean that the evaluating the error estimator has substantially higher costs that computing the problem itself. High order reconstructions of low order solutions usually give optimal results, see~\cite{RichterWick2015,BeckerRannacher2001}	
	
	It is essential that the adjoint solution is approximated in a space that is somehow orthogonal to the primal one, but this approximation must not necessarily be of higher accuracy.
\end{remark}
In what follows, we employ the DWR method to construct a functional error estimator depending solely on the neural network approximation of Poisson equation. 
\subsection{Dual weighted residual estimation for PINNs}

The error identity~\eqref{errorrepresentation} only requires that the discrete approximation is $H^1$-conforming, i.e. that $\tilde u\in \tilde V\subset H^1(\Omega)$ comes from a subspace. We have however in no way used that $\tilde u$ actually solves the problem. This is also the basis for incorporating iteration errors into the estimator, see~\cite{Meidner2009,Rannacher2013}. It further means that the error representation is directly applicable to any type of approximation method, e.g., a finite difference approximation can be considered as interpolation into a finite element space, or, in the context of this work, identity~\eqref{errorrepresentation} can directly be used to estimate errors from PINN representations as done in~\cite{MR}. 

However, in~\cite{MR} the adjoint solution $z$ is approximated with a very coarse finite element approximation. The surprising observation from~\cite{MR} is that extremely coarse spaces $z_H\in V_H\subset H^1(\Omega)$ consisting of only very few basis functions were sufficient to estimate the functional error $J(u-u_{\lambda,{\cal N}})$ during the complete training process. However, the use of a (albeit coarse) finite element mesh to approximate the dual solution negates any potential advantage of PINNs, for example in the treatment of very high-dimensional problems.
In~\cite{RothSchroederWick2022} the opposite path was taken, where the primal solution was finite element based and the adjoint solution approximated by a network. The practical application of this approach however suffers from the same inadequacy as our attempt, as the worst of both worlds is combined, the inflexibility of mesh-based finite elements with the bad approximation property of PINNs.

In the following, we assume that the adjoint solution itself is approximated by a neural network. We consider the energy formulation of the Poisson equation~\eqref{loss:deepritz} and introduce the adjoint functional
\[
E_{\lambda,\text{adj}}(z_\theta) = \frac{1}{2}\int_\Omega|\nabla z_\theta(x)|^2\,\text{d}x - \int_\Omega j(x)z_\theta(x)\,\text{d}x + \frac{\lambda}{2}\|z_\theta\|^2_{\partial\Omega}
\]
which, if minimized over $H^1(\Omega)$, corresponds to the variational problem
\[
z_\lambda\in H^1(\Omega):\quad
\langle \nabla \phi,\nabla z_\lambda\rangle_\Omega 
+\lambda \langle \phi,z_\lambda\rangle_{\Omega} = \langle j,\phi\rangle_{\Omega}=J(\phi)\quad\forall \phi\in H^1(\Omega).
\]
For the error between $u_\lambda\in H^1(\Omega)$ given by~\eqref{laplace1:dist} and any type of neural network approximation $u_{\lambda,\mathcal{N}}\in V_{\cal N}$ this gives

\begin{equation}\label{errorrepresentation:NN}
\begin{split}
J(u_\lambda-u_{\lambda,\cal N})  & = 
\langle \nabla (u_\lambda-u_{\lambda,\mathcal{N}}),\nabla z_\lambda\rangle_\Omega 
+\lambda \langle u_\lambda-u_{\lambda,\mathcal{N}},z_\lambda\rangle_{\partial\Omega}\\
& = \langle f,z_\lambda\rangle_\Omega - \langle \nabla u_{\lambda,\mathcal{N}},\nabla z_\lambda\rangle_\Omega 
-\lambda \langle u_{\lambda,\mathcal{N}},z_\lambda\rangle_{\partial\Omega}
\end{split}
\end{equation}
A first computable error approximation is obtained by simple replacing the unknown adjoint $z_\lambda\in H^1(\Omega)$ by its neural network approximation $z_{\lambda,\mathcal{N}}$
\begin{equation}\label{est1}
\eta(u_{\lambda,\mathcal{N}},z_{\lambda,\mathcal{N}})=\langle f,z_{\lambda,\mathcal{N}}\rangle_\Omega - \langle \nabla u_{\lambda,\mathcal{N}},\nabla z_{\lambda,\mathcal{N}}\rangle_\Omega 
-\lambda \langle u_{\lambda,\mathcal{N}},z_{\lambda,\mathcal{N}}\rangle_{\partial\Omega}.
\end{equation}
This estimator is highly heuristic as replacing $z_\lambda$ by $z_{\lambda,\mathcal{N}}$ introduces an unknown error. 

\begin{remark}[Choice of the adjoint network]
So far, we did not comment on the network architecture for the adjoint approximation $z_{\lambda,\mathcal{N}}$. If we draw on finite element theory, it is crucial that $z_{\lambda,\mathcal{N}}\not\in V_h$, i.e. that the dual approximation does not lie in the primal test space. From the discussion in Section~\ref{subsection2.2} however we noted that the primal approximation $u_{\lambda,\mathcal{N}}\in V_{\cal N}$ can be considered as variational solution, however, tested with functions coming from a distinct test space $V'_{\cal N}u_\theta)\neq V_{\cal N}$. Thereby, we will approximate $z_{\lambda,\mathcal{N}}$ in~\eqref{est1} within the same network set $V_{\cal N}$. It is even possible to train both $u_{\lambda,\mathcal{N}}$ and $z_{\lambda,\mathcal{N}}$ within the same network at the same time. In Section~\ref{section5} we demonstrate that this very simple estimator gives highly accurate results  from early in the training process.

\end{remark}

 Since the estimator obtained in the previous section lacks pointwise accuracy, it is necessary to improve it. To this end, the next two sections successively present the localization technique commonly used in the finite element literature to enhance the estimator, and then adapt this approach to construct an improved estimator for PINNs and the Deep Ritz method based solely on neural network approximations.

\section{Adaptivity}\label{section4}

In this section we describe the use of the error estimator for steering the training process in an adaptive way. Comparable to $h$-adaptivity in the finite element method we aim at a localization of the estimator~\eqref{est1} which is then used to sample quadrature nodes for the approximation of the loss term.

\subsection{Localizations in the Finite Element case}
The key point in localizing an error estimate as (here written in the standard finite element setting)
\begin{equation}\label{est:fe}
\eta(u_h,z) = \langle f,z\rangle_{\Omega} - \langle\nabla u_h,\nabla z\rangle_\Omega - \lambda\langle u_h,z\rangle_{\partial\Omega}
\end{equation}
is a splitting into positive local quantities $\eta_i$ which we call \emph{error indicators}
\[
\eta(u_h,z) \approx \sum_i \eta_i.
\]
These $\eta_i$ can be element- or node-based and are then used for adaptivity. The local quantities are called \emph{error indicators}. Since the error estimator $\eta(u_h,z)$ is usually oscillating and changing its sign, the main goal of localization is that the \emph{indicator index} 
\begin{equation}\label{ind}
I_\text{ind}:=\frac{\sum_i|\eta_i|}{|\eta(u_h,z)|} \approx 1
\end{equation}
is bounded as close to one as possible. In~\cite{RichterWick2015} we have discussed several approaches to accomplish this goal. The first step is usually to apply Galerkin-Orthogonality and rewrite~\eqref{est:fe} as
\begin{equation}\label{est:fe1}
\eta(u_h,z) = \langle f,z-I_h z\rangle_{\Omega} - \langle\nabla u_h,\nabla (z-I_h z)\rangle_\Omega - \lambda\langle u_h,z-I_h z\rangle_{\partial\Omega}
\end{equation}
where $I_h:H^1(\Omega)\to V_h$ is an interpolation into the finite element space. To get a computable quantity, we next approximate the adjoint $z\in H^1(\Omega)$, for instance by reconstruction of the discrete adjoint $z_h\in V_h$ in a higher accuracy space $z_h^{(*)}=\pi_h^{(*)}z_h\in V_h^{(*)}$
\begin{equation}\label{est:fe:approx1}
\eta(u_h,z)\approx \eta(u_h,z_h^{(*)})
 =
 \langle f,\pi_h^{(*)}z_h-z_h\rangle_{\Omega} 
 - \langle\nabla u_h,\nabla (\pi_h^{(*)}z_h-z_h)\rangle_\Omega 
 - \lambda\langle u_h,\pi_h^{(*)}z_h-z_h\rangle_{\partial\Omega}
\end{equation}

In~\cite{RichterWick2015} it is discussed that, in the finite element case, this is not sufficient to reach $I_\text{ind} < C$ independent of $h$. For example, the obvious choice of a localization towards the elements $T\in\Omega_h$ of the mesh
\begin{equation*}
\eta_T:= \langle f,\pi_h^{(*)}z_h-z_h\rangle_{T} 
 - \langle\nabla u_h,\nabla (\pi_h^{(*)}z_h-z_h)\rangle_T 
 - \lambda\langle u_h,\pi_h^{(*)}z_h-z_h\rangle_{\partial\Omega\cap \partial T}
\end{equation*}
results in $I_{\text{ind} }= O(h^{-1})$. Remedy is given by transferring the strong residual and introducing jump terms of the normal derivative, see~\cite{BeckerRannacher1995} for this standard approach. A computationally simple and highly accurate approach consists of introducing a partition of unity $\sum_i\psi_i\equiv 1$ and defining the local indicators as
\begin{equation}\label{est:fe:approx}
\eta_i:=
 \langle f,(\pi_h^{(*)}z_h-z_h)\psi_i\rangle_{\Omega} 
 - \langle\nabla u_h,\nabla \big((\pi_h^{(*)}z_h-z_h)\psi_i\big)\rangle_\Omega 
 - \lambda\langle u_h,(\pi_h^{(*)}z_h-z_h)\psi_i\rangle_{\partial\Omega}.
\end{equation}

\subsection{Localization of the PINN (DRM)  estimator}

The essential step in localization of the finite element estimator is the introduction of error terms $z-\psi_h$, usually $\psi_h = I_h z$. Classically, these factors are called the \emph{dual weights} of the estimator~\cite{BeckerRannacher1995}. This idea is transferred to the Deep Ritz method by using the special structure of the variational setting described in Lemma~\ref{lemma:var}. 

Assuming that $u_{\lambda,\mathcal{N}}\in V_{\cal N}$ and $z_{\lambda,\mathcal{N}}\in V_{\cal N}$ are neural network approximations to primal and adjoint solution, the heuristic error estimator~\eqref{est1} will be localized by using an orthogonality relation in the spirit of~\eqref{est:fe}. If $u_{\lambda,\mathcal{N}}$ would be the exact minimizer of the energy functional~\eqref{loss:deepritz} in $V_{\cal N}$ (meaning that we neglect the optimization and the quadrature error), it would be characterized by the Euler-Lagrange equation~\eqref{Euler-Lagrange}, which defines the test set as 
\[
    V_\mathcal{N}^\prime(u_{\lambda,\mathcal{N}}) := \left\{ x\mapsto  \frac{d}{ds} u_{\theta+s\beta}(x)\big|_{s=0}: \quad  u_\theta= u_{\lambda,\mathcal{N}}, \quad \beta \in \Theta=\mathbb{R}^{n_p}\right\}.
\]
Hence, for any $z_\omega'\in V_{\cal N}'(u_{\lambda,\mathcal{
N}})$ it holds
\begin{equation}\label{est:final}
\begin{split}
\eta(u_{\lambda,\mathcal{N}},z_{\lambda,\mathcal{N}}) & = 	\eta(u_{\lambda,\mathcal{N}},z_{\lambda,\mathcal{N}},z'_{\lambda,\omega}) \\
& =
	\langle f,z_{\lambda,\beta}-z'_{\lambda,\omega}\rangle_\Omega - \langle \nabla u_{\lambda,\mathcal{N}},\nabla (z_{\lambda,\beta} - z_{\lambda,\omega}')\rangle_\Omega 
-\lambda \langle u_{\lambda,\mathcal{N}},z_{\lambda,\mathcal{N}}- z'_{\omega}\rangle_{\partial\Omega}.
\end{split}
\end{equation}
This identity~\eqref{est:final} corresponds to the variational formulation of the Poisson problem. In Lemma~\ref{lemma:var} we derived exactly this equation for defining the minimizer of the Energy form and used it in connection with the Deep Ritz method. So far,~\eqref{est:final} is agnostic with respect to the specific scheme for training the neural network. We can consider both the Deep Ritz method and the classical PINN as an approximating scheme for the solution~\eqref{Euler-Lagrange}.

In applications, ``identity''~\eqref{est:final} will not he exactly satisfied due to multiple reasons. First of all, the Deep Ritz method must approximate all integrals with numerical quadrature, which introduced a substantial error. Second, we usually are not able to find a good minimum. Furthermore, we intend to use the error estimator during training to steer the learning process. This means that the error estimator is applied early in training, when $u_\theta$ might still be far away from the optimum. 

\eqref{est:final} is the basis for our localization and we introduce the \emph{local indicator function}
\begin{equation}\label{est:local}
\mu(x):=\begin{cases}
  f(x)\cdot \big(z_{\lambda,\beta}(x)-z'_{\lambda,\omega}(x)\big)\\
  \qquad 
  -  \nabla u_{\lambda,\mathcal{N}}(x)\cdot \nabla \big(z_{\lambda,\beta}(x)
  - z_{\lambda,\omega}'(x)\big) & x\in\Omega\\
  -\lambda u_{\lambda,\mathcal{N}}(x)\cdot\big(z_{\lambda,\mathcal{N}}(x)- z'_{\omega}(x)\big) & x\in\partial\Omega
\end{cases} 
\end{equation}

We will see with numerical test cases in Section \ref{section5} that this ``first step of localization'' is sufficiently accurate to deliver both an estimate and information on the localization. 

In the finite element context, localization by use of Galerkin orthogonality is not sufficient, compare~\cite{RichterWick2015} and a second step, e.g. using the partition of unity localization~\cite{RichterWick2015} is required. This could easily be applied to~\eqref{est:final} as well. However, we will observe that the estimator and its localization is already suffiicently accurate.

For comparison, numerical examples will also be run for localizations based on the simple estimator
\begin{equation}\label{est:local:simple}
\mu_\text{simple}(x):=\begin{cases}
  f(x)\cdot z_{\lambda,\beta}(x) -  \nabla u_{\lambda,\mathcal{N}}(x)\cdot \nabla z_{\lambda,\beta}(x) & x\in\Omega\\
  -\lambda u_{\lambda,\mathcal{N}}(x)\cdot z_{\lambda,\mathcal{N}}(x) & x\in\partial\Omega
\end{cases} 
\end{equation}

To obtain this estimator, the property that the energy minimizer found using the Deep-Ritz method is the one used; this relationship derived from the energy allows us to perform localization in order to obtain this estimator.   We can also consider this estimator for Pinns problems, since an approximation of a Laplace problem using the Pinns method is also an approximation via the Deep Ritz problem (approximation of the same problem but using different methods); therefore, this approximation via PINNs is also a minimizer of the energy, even if this minimizer is local.  

\subsection{Description of adaptive method}\label{subsection-description}

In this section, we describe how to numerically compute the estimator derived in the previous section and demonstrate how it can be used for adaptivity. First, we compute the estimator numerically. Next, we assess its accuracy and precision. Finally, we use it for adaptivity. 
  The computation of the estimator is divided into three steps, which are described as follows.
  \begin{enumerate}
    \item We train a neural network to learn the adjoint equation. At this stage, we obtain a PINN/Deep Ritz approximation of the solution to the adjoint equation. Then we obtain  $z_{\lambda, \mathcal{N}}$.  
      \item We start training a neural network $u_\theta$  to learn the primal equation.\\
       This step is important because the estimator can only be computed during this training phase. In practice, we train the neural network for a few epochs so that it can start learning the equation. The training points used at this stage are sampled randomly, as in the classical approach. So the preliminary solution we obtain is denoted by $u_{\lambda,\theta}$. 
      \item  Consider a neural network $z^\prime_\omega$. Freeze the parameters of the first layer of  the neural network $z^\prime_\omega$ and fix them equal to the parameters of the neural network $u_{\lambda,\theta}$ as shown in Example \ref{Example}. $z^\prime_\omega \in V_\mathcal{N}^\prime  (u_{\lambda,\theta})$.   And finally, we train the neural network $z^\prime_\omega$ with the freezing parameters to learn the adjoint equation. We obtain the solution $z'_{\omega}$. 
  \end{enumerate}
  Therefore, using the preliminary solution $u_{\lambda,\theta}$,  and the different solutions of the adjoint equation $z_{\lambda,\mathcal{N}}$ and $z^\prime_\omega$, we  compute the estimator defined in \eqref{est:final}. The complete procedure is summarized in Algorithm~\ref{alg:adaptive}.

\begin{algorithm}
\begin{algorithmic}[1]
\Require Batch size $N$; stopping criterion $\mathcal{S}$
\Ensure Trained primal $u_\theta$ and adjoint $z_{\theta^\prime }$

\State $X \gets \textsc{SampleUniform}(N)$
\State $z_{\lambda,\mathcal{N}} \gets \textsc{InitTrainAdjoint}(z_{\theta^\prime} )$
\State $u_\theta \gets \textsc{InitPrimal}()$

\While{not $\mathcal{S}(u_\theta, z_{\lambda,\mathcal{N}})$}
    \State $u_{\lambda,\theta} \gets \textsc{TrainPrimal}(u_\theta, X)$
    \State $z_\omega^\prime \gets \textsc{InitAdjoint}(u_\theta)$
    \State $z_\omega^\prime \gets \textsc{TrainAdjoint}(z_\omega^\prime)$
    \State $E \gets \textsc{EstimateError}(u_{\lambda,\theta}, z_{\lambda,\mathcal{N}}, z_\omega^\prime)$
    \State $X \gets \textsc{SampleFromError}(E, N)$
    % \Comment{e.g., importance resampling / refinement}
    \State $u_\theta \gets u_{\lambda, \theta} $  
    \State  $u_{\lambda,\mathcal{N}} \gets \textsc{TrainPrimal}(u_\theta, X)$ 
\EndWhile

\State \Return $(u_{\lambda,\mathcal{N}})$
\end{algorithmic}
\caption{Training with Adaptive Resampling}
\label{alg:adaptive}
\end{algorithm}

\paragraph{Adaptive sampling method}
We consider the measure 
\[
p_\mu : x \mapsto \frac{|\mu(x)|)}{\int_\Omega|\mu (x)| }.
\]
The adaptive sampling consist of sampling using the measure $p_\mu$. This sampling depends only on the error estimator $\mu$, and it helps to sample more points with large functional error. 
  After sampling new points with our probability measure, we do not change the loss even though we replace all the collocation points. So for the rest of the training; the loss function of PINNS and Deep-Ritz Method are respectively:
  
\begin{equation}
E_{\lambda, M, N} := \frac{1}{N}\sum_{i=1}^{N} \left| -\Delta u_\theta (x_i) - f(x_i) \right| ^2
+  \frac{\lambda}{M} \sum_{j=1}^{M}\left|u_\theta(y_j)-g(y_j)\right|^2, 
\end{equation}
where $x_i  \sim p_\mu $ and $y_j $ is uniformly sampled and
\begin{equation}\label{DR:loss-adap}
    E_{\lambda, M, N}(u_\theta) := \frac{\vert \Omega \rvert }{N}\sum_{i=1}^{N} \frac{1}{2} \lvert \nabla u_\theta (x_i) \rvert ^2 - f(x_i) u_\theta(x_i)
    + \lambda \frac{\vert \partial \Omega \rvert }{2M} \sum_{j=1}^{M}|u_\theta(y_j)-g(y_j)|^2 ,   
\end{equation}
with $x_i \sim p_\mu $ and $y_j  $ uniformly sampled.

A second approach is to add new collocation points to the collocations used during the pre-training phase. This approach is also called adaptive refinement. The points to be added are obtained by sampling with respect to the new distribution $
p_\mu$. This implies that we only add more points where the functional error is large.
%  \begin{comment}
 %    {\color{orange}Some more details in this paragraph. Also exact formular what is being used point-wise. }
 % Instead of selecting the training data randomly, we have at this point a new measure that helps to select the training data (points) where the functional error is large.
 % \end{comment}
This method  guarantees to find a PINN / Deep-Ritz approximation with a smaller (optimal) value of the functional $J$. We illustrate this in the following section with some examples. 

 \section{Numerical experiments}\label{section5}

In the following we present numerical examples solving the Poisson problem $-\Delta u=f$ in $\Omega$ with $u=g$ on $\partial\Omega$ for several low- and high-dimensional configurations. The focus of the presentation is on the accuracy of the error estimator to be used as possible stopping criterion and, in particular, on the adaptive sampling of collocation points. We will test both the standard PINN based on the strong formulation and the Deep Ritz method.

\subsection{Application to the PINNs framework}

In the first examples we will consider the PINN approach based on minimizing the classical residual of the equation. As examples we consider the Poisson problem in 2 and in higher dimensions. 

\subsubsection{Case 1: Low Dimensional Poisson Problem with smooth solution}

On the domain $\Omega=(-1,1)^2\subset\mathbb{R}^2$ we consider the problem
\[
-\Delta u = \frac{\pi^2}{2}
\cos\left( \frac{\pi}{2} x \right)
\cos\left( \frac{\pi}{2} y \right)\text{ in }\Omega,\quad u=0\text{ on }\partial\Omega.
\]
The analytical solution is given by
\[
u(x,y) = 
\cos\left( \frac{\pi}{2} x \right)
\cos\left( \frac{\pi}{2} y \right).
\]
As functional of interest we choose the average of the solution over the  complete domain, i.e.
\[
 J(u) = \int_\Omega u(x)\,\text{d}x
\]  
which gives rise to the adjoint problem
\[
-\Delta z =1\text{ in }\Omega,\quad z=0\text{ on }\partial\Omega.
\]
We use a fully residual neural network composed of 2 blocks, each consisting of two fully connected layers, with 64 neurons per layer, to approximate the exact solution of the equations. The neural networks are trained for 2000 epochs. 

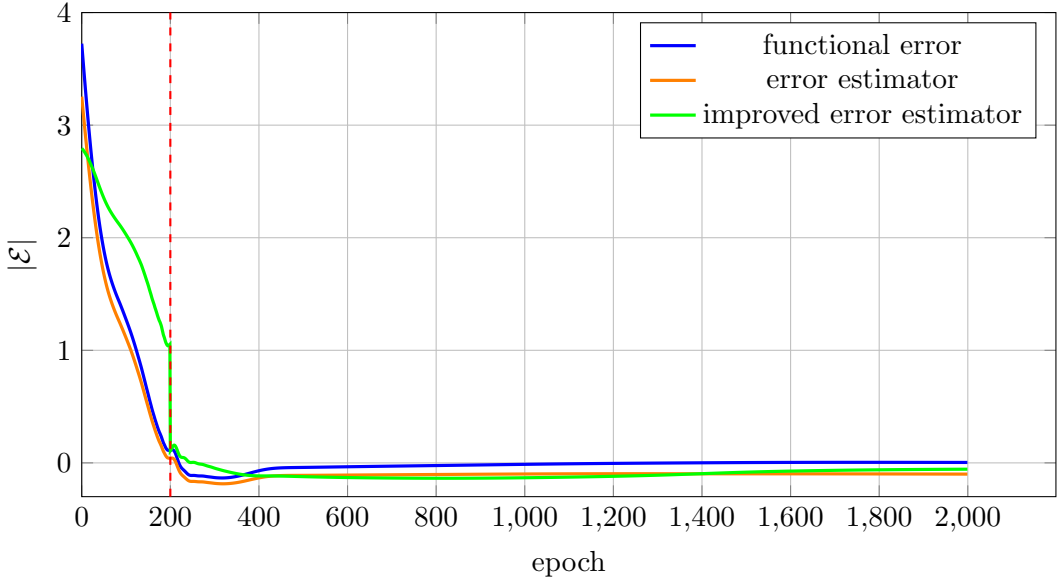
\begin{figure*}
\centering
\begin{tikzpicture}
\begin{axis}[
  width=0.9\linewidth,
  height=8cm,
  xlabel={epoch},
  ylabel={$|{\cal E}|$},          % souvent utile pour une loss
  grid=both,
  xmin=0,
  ymin=-0.3,
  ymax=4.0
]
\addplot[blue,very thick]
table[
  col sep=comma,
  x=epoch,
  y=Error,
]{adaptive_pinn_ex2.csv};
\addlegendentry{functional error}
\addplot[orange,very thick] table[
  col sep=comma,
  x=epoch,
  y =Est Error,  mark=none
]{adaptive_pinn_ex2.csv};
\addlegendentry{error estimator}
\addplot[green,very thick] table[
  col sep=comma,
  x=epoch,
  y =Est Error prime,
  mark=none
]{adaptive_pinn_ex2.csv};
\addlegendentry{improved error estimator}
% \node[blue] at (axis cs:900,1e-4) {Standard method};
\addplot[red, thick, dashed]
  coordinates {(200,\pgfkeysvalueof{/pgfplots/ymin})
               (200,\pgfkeysvalueof{/pgfplots/ymax})};
\end{axis}
\end{tikzpicture}
\caption{Case 1: Evolution of the functional error and its estimators during training. The vertical dashed red line indicates the iteration at which resampling is performed. }
\label{fig:figure1}
\end{figure*}

Figure~\ref{fig:figure1} shows the functional error $J(u)-J(u_\theta)$, the error estimator \eqref{est1} as well as the localized estimator~\eqref{est:final} during training. We observe that the neural network estimator~\eqref{est1} is highly accuracy already very early during training. The DWR method does not deliver a strict upper bound but an error approximation that also includes the sign of the error. Fig.~\ref{fig:figure1} shows a change of sign, approximately, at epoch 220 and the estimator follows this behavior accurately. The localized estimator~\eqref{est:final} (in green) is less accurate early in training. However, it is still close to the real error and gives superior results in the later states of training.
 
 \begin{figure*}\label{pointwise:est}
    \centering
    \includegraphics[width=1.0\textwidth]{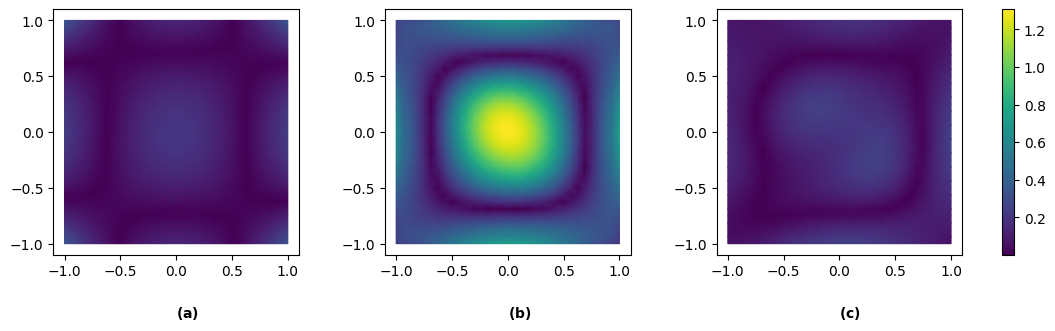}
   \caption{
(a) Pointwise functional error, (b) corresponding pointwise error estimator (b), (c) corresponding pointwise improved error estimator with localization. }
    \label{fig:figure2}
\end{figure*}

Figures \ref{fig:figure2}\,(a) shows the error $u(x)-u_\theta(x)$ after training of the neural network. For comparison, Fig.~\eqref{fig:figure2}\,(b) shows the local values of the error estimator $\mu_\text{simple}(x)$, see~\eqref{est:local:simple}. Although the overall estimate is highly accurate, i.e.
\[
\frac{\int_\Omega u(x)-u_\theta(x)\,\text{d}x}{\int_\Omega \mu_\text{simple}(x)\,\text{d}x}\approx 1,
\]
the localization of the error shows excessive variations in the positive and negative directions. Used for adaptivity, which means steering the adaptive sampling in the context of PINNs, a faulty selection of collocation points might result. If we consider the localized estimator $\mu(x)$ shown in Fig.~\ref{fig:figure2}\,(c), we observe much smaller amplitudes in the error while the overall efficiency is preserved, i.e. it still holds
\[
\frac{\int_\Omega u(x)-u_\theta(x)\,\text{d}x}{\int_\Omega \mu(x)\,\text{d}x}\approx 1,
\]
as can be seen in Figure~\ref{fig:figure1}.

\begin{remark}[Error localization]
A direct comparison of the local error $u(x)-u_\theta(x)$ in Fig.~\ref{fig:figure2}\,(a) and the error estimators $\eta(x)$ and $\eta_\text{simple}(x)$ in Fig.~\ref{fig:figure2}\,(b) and (c) must be made with caution. The DWR estimator does not provide a distribution of the local error but rather indications of the origin of the error. The local indicators should be regarded as sensitivities of the total error, which in particular take into account the non-locality of the finite element approximation due to the pollution effect.  
\end{remark}

\begin{figure*}
    \centering
    \includegraphics[width=0.8\textwidth]{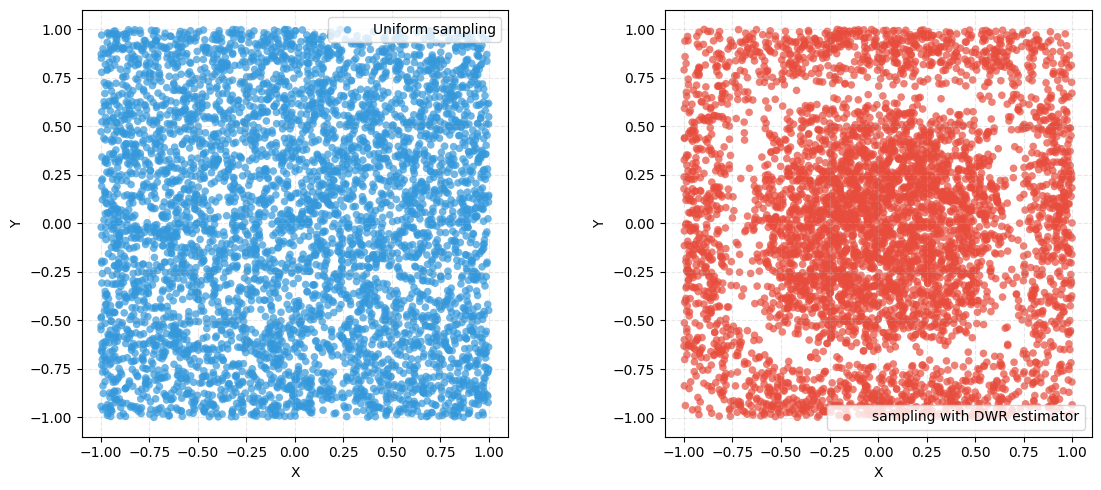}
    \caption{Case 1: Comparison of uniform sampling (left) vs. sampling using the DWR measure $p_\mu$ (right).}
    \label{fig:Sampling}
\end{figure*}

\begin{comment}
\begin{figure}[H]
    \centering
    \includegraphics[width=0.8\textwidth]{Estimation after localisation.png}
    \caption{Evolution of the functional error and the improved estimator during training.}

    \label{fig:figure3}
\end{figure}

\begin{figure}[H]
    \centering
    \includegraphics[width=0.8\textwidth]{Exo2_est_prime.png}
   \caption{Pointwise functional error and the improved estimator.
Pointwise functional error (a) and corresponding improved pointwise error estimator (b) 
over the computational domain $(x,y)$. }

    \label{fig:figure4}
\end{figure}
    
\end{comment}
After computing the improved error estimator at epoch 200, the training points are resampled based on this estimator and used to continue the training of the neural network.
Figure~\ref{fig:Sampling} shows the difference between the uniform distribution of collocation points and the one generated by means of the DWR estimator. Both cases involve about 5000 points. The DWR measure $p_\mu$ based  sampling generates a distribution which is aligned with regions of larger values of the indicator $|\mu(x)|$, compare Fig.~\ref{fig:Sampling}.

\begin{figure*}
\centering
\begin{tikzpicture}
\begin{axis}[
  width=0.9\linewidth,
  height=8cm,
  xlabel={epoch},
  ylabel={$|{\cal E}|$},
  ymode=log,           % souvent utile pour une loss
  grid=both,
  xmin=0,
  ymin=2.1e-6,
  ymax=10
]
\addplot[blue]
table[
  col sep=comma,
  x=epoch,
  y expr={abs(\thisrow{Error})},
]{adaptive_pinn_func_plot_ex2.csv};
\addlegendentry{Adaptive method}
\addplot[orange] table[
  col sep=comma,
  x=epoch,
  y expr={abs(\thisrow{Functional error})},
  mark=none
]{result_direct_pinn_ex2.csv};
% \node[blue] at (axis cs:900,1e-4) {Standard method};
\addplot[red, thick, dashed]
  coordinates {(200,\pgfkeysvalueof{/pgfplots/ymin})
               (200,\pgfkeysvalueof{/pgfplots/ymax})};
\addlegendentry{Standard method}
\end{axis}
\end{tikzpicture}
\caption{Case 1: Absolute value of the functional error 
$|\mathcal{E}(u)| = |J(u) - J(u_\theta)|$ during the training. The vertical dashed red line indicates the iteration at which resampling is performed.}
\label{fig:jduringtraining}
\end{figure*}
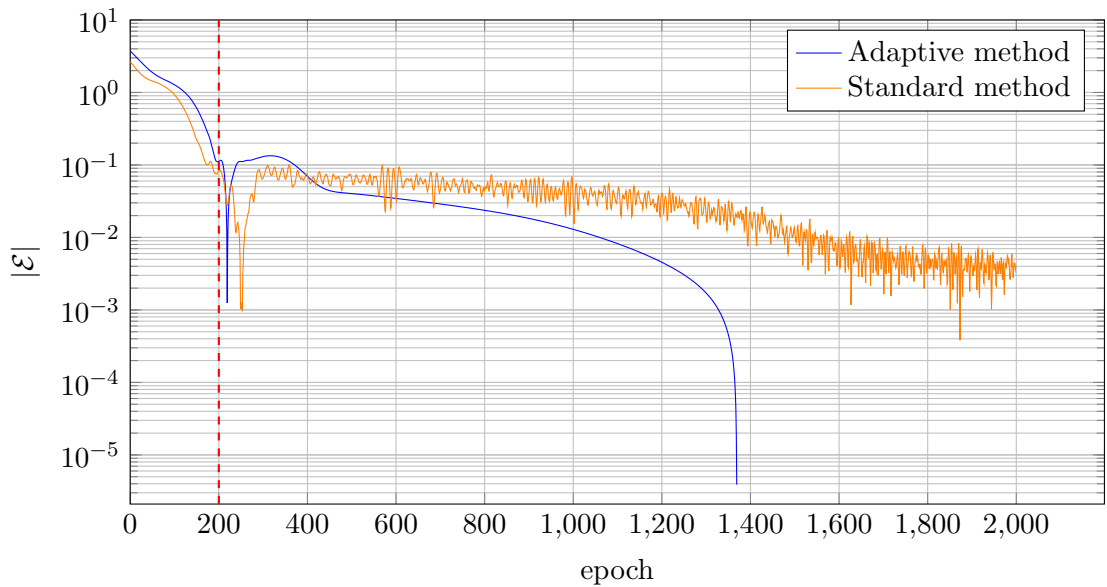

Finally, figure*~\ref{fig:jduringtraining} the evolution of the functional error $|J(u)-J(u_\theta)|$ during the learning process. After adaptive sampling using the DWR estimator, the functional error converges faster and reaches a smaller value than in the non-adaptive training. The discrepancy between te functional error at the beginning comes from the usual randomness of machine learning approaches originating, e.g., in the random initialization of the parameters.

%\end{Case}

\subsubsection{Case 2: Low dimensional problem with non-smooth functional}
We consider a similar primal problem given by 
\begin{equation}
-\Delta u= 2 \sin(x) \sin(y)\text{ in }\Omega=(0,\pi)^2,\quad u=0\text{ on }\partial\Omega
\end{equation}
and the average of the solution in a subset as error functional:
\[
 J(\phi)=\int_C \phi(x)\,\text{d}x,\quad C:=\{(x,y)\in \mathbb{R}^2:\;
 (x-\frac{\pi}{2})^2+(y-\frac{\pi}{2})^2\le 1\}.
\] 
This gives rise to the adjoint problem
\[
-\Delta z = 1_\mathcal{C}\text{ in }\Omega,\quad z=0\text{ on }\partial\Omega,
\]
which has some reduced regularity along the interface between $C$ and $\Omega\setminus C$. 

We use a fully residual neural network composed of four blocks, each consisting of two fully connected layers, with 32 neurons per layer, to approximate the exact solution of the equations. The neural networks are trained for 5000 epochs.

\begin{figure*}
\centering
\begin{tikzpicture}
\begin{axis}[
  width=0.9\linewidth,
  height=8cm,
  xlabel={epoch},
  ylabel={$|{\cal E}|$},
  ymode=log,           % souvent utile pour une loss
  grid=both,
  xmin=0,
  ymin=4.e-3,
  ymax=10
]
\addplot[blue]
table[
  col sep=comma,
  x=epoch,
  y expr={abs(\thisrow{Error})},
]{adaptive_pinn_func_plot_ex1.csv};
\addlegendentry{Adaptive method}
\addplot[orange] table[
  col sep=comma,
  x=epoch,
  y expr={abs(\thisrow{Functional error_standard})},
  mark=none
]{adaptive_pinn_func_plot_ex1.csv};
% \node[blue] at (axis cs:900,1e-4) {Standard method};
\addplot[red, thick, dashed]
  coordinates {(300,\pgfkeysvalueof{/pgfplots/ymin})
               (300,\pgfkeysvalueof{/pgfplots/ymax})};
\addlegendentry{Standard method}
\end{axis}
\end{tikzpicture}
\caption{Case 2: Absolute value of the functional error 
$|\mathcal{E}(u)| = |J(u) - J(u_{\mathcal{N}})|$ during the training. The vertical dashed red line indicates the iteration at which resampling is performed.}
\label{fig:case2}
\end{figure*}
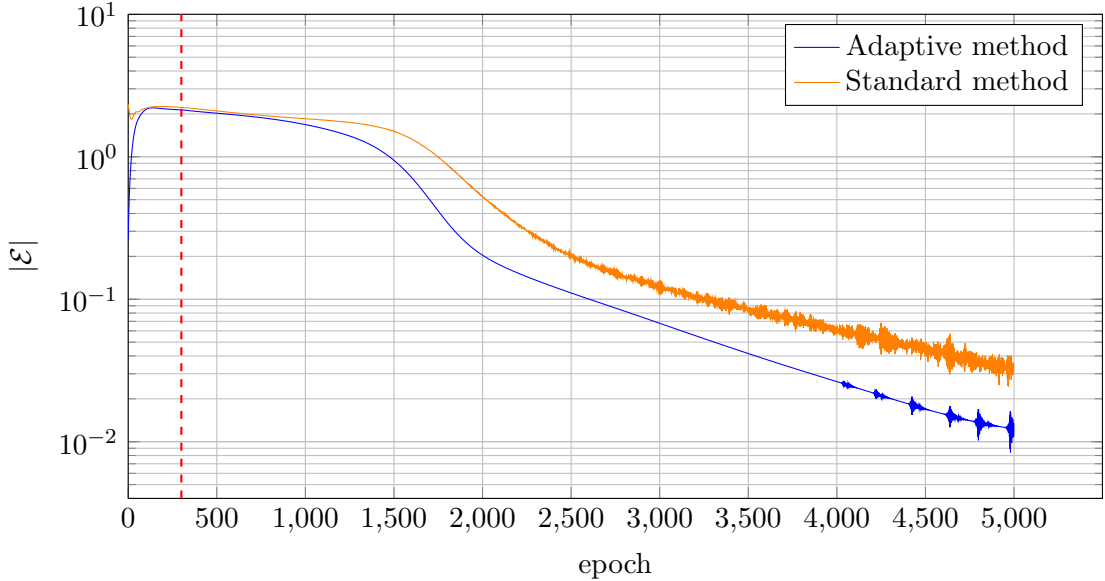

For the adaptive method, the training is first performed using randomly selected collocation points, replaced by adaptive sampled points based on the estimator at epoch 300.  The results shown in Figure~\ref{fig:case2} again demonstrate the improvement of the error reduction rate when adaptive sampling is applied. Also, the error level reached at the end of training is substantially improved. 

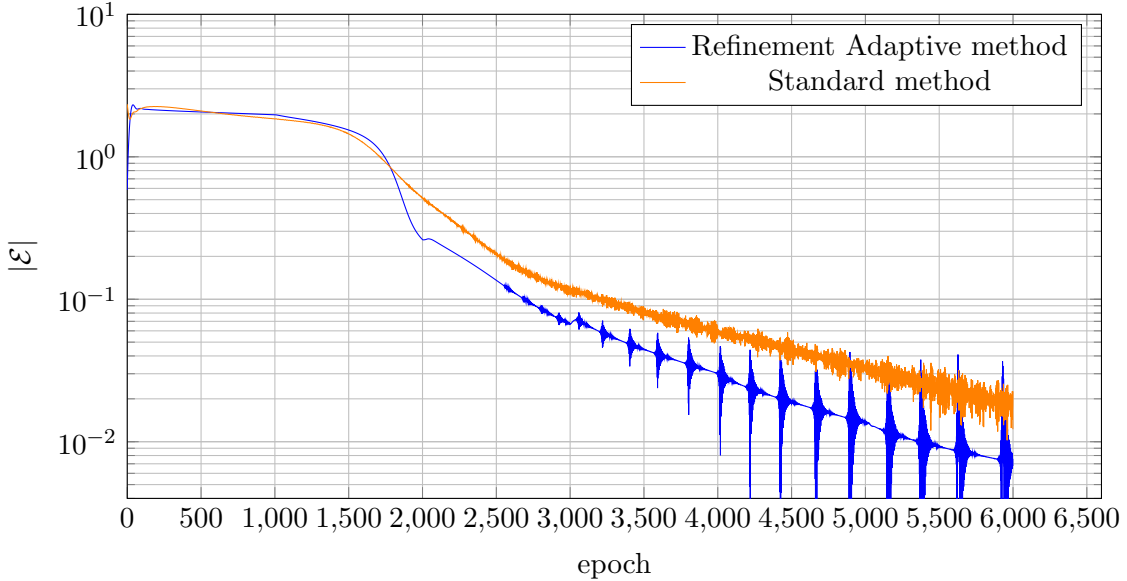
\begin{figure*}
\centering
\begin{tikzpicture}
\begin{axis}[
  width=0.9\linewidth,
  height=8cm,
  xlabel={epoch},
  ylabel={$|{\cal E}|$},
  ymode=log,           % souvent utile pour une loss
  grid=both,
  xmin=0,
  ymin=4.e-3,
  ymax=10
]
\addplot[blue]
table[
  col sep=comma,
  x=epoch,
  y expr={abs(\thisrow{Error_adaptive})},
]{DWR-D-adaptive_pinn_ex1.csv};
\addlegendentry{Refinement Adaptive method}
\addplot[orange] table[
  col sep=comma,
  x=epoch,
  y expr={abs(\thisrow{Error-Non_Adaptive})},
  mark=none
]{DWR-D-adaptive_pinn_ex1.csv};
% \node[blue] at (axis cs:900,1e-4) {Standard method};
\addlegendentry{Standard method}
\end{axis}
\end{tikzpicture}
\caption{Case 2: Comparison of the standard method with successive adaptive sampling after 100 epochs each. Absolute value of the functional error 
$|\mathcal{E}(u)| = |J(u) - J(u_{\mathcal{N}})|$ during the training.}
\label{fig:case2-refinement}
\end{figure*}

Let us still consider the same Poisson equation and try to apply the adaptive sampling method. We use the same neural network architecture. We start the training with 200 collocation points and every 1000 epochs we sample 500 points according to the probability distribution $p_\mu$ that we add to the training points. 
Figure~\ref{fig:case2-refinement} shows the evolution of the functional error as a function of the number of epochs.

\newpage 

\subsubsection{Case 3: High-dimensional problem}
We consider the Poisson problem $-\Delta u=f$ in the 5 dimensional annulus
\[
\Omega = B_2(0)\setminus B_1(0),\quad B_r(0):=\{x\in\mathbb{R}^5\;:\; \|x\|_2\le r\}
\]
with homogeneous boundary values $u=0$ on $\partial\Omega$.
For the analytical solution 
\[
u(x) = \frac{\tanh\big( \gamma(x)\big)}{800 \tanh(10)} -\frac{\gamma(x)}{80},\quad
\gamma(x):=20\left(\|x\|_2-\frac{3}{2}\right),
\]
the right hand side $f=-\Delta u$ is given as
\begin{equation*}
 f(x) = -\frac{1}{800}
\left\{
-\frac{800}{\tanh(10)}\,
\frac{\tanh\big(\gamma(x)\big)}{\cosh^2\big(\gamma(x)\big)}
+\frac{4}{\|x\|_2}
\left(
\frac{20}{\tanh(10)}\cdot 
\frac{1}{\cosh^2\big(\gamma(x)\big)}-2
\right)
\right\}.
\end{equation*}
We consider the error functional
\begin{equation*}
   J(u) = \int_\Omega f(x)u(x) dx .
\end{equation*} 
which is a weighted average of the solution. For this special choice, the adjoint problem gets $-\Delta z=f$ with $z=u$, which eases the validation.

We use a fully connected network composed of 6 layers, each consisting of 32 neurons per layer, to approximate the exact solution of the equations. The neural networks are trained for 5000 epochs. For the adaptive method, the training is first performed using randomly selected points. At epoch 400, we compute the functional error estimator, we resample the interior points, and we use them to continue the training of the neural network to approximate the primal equation.

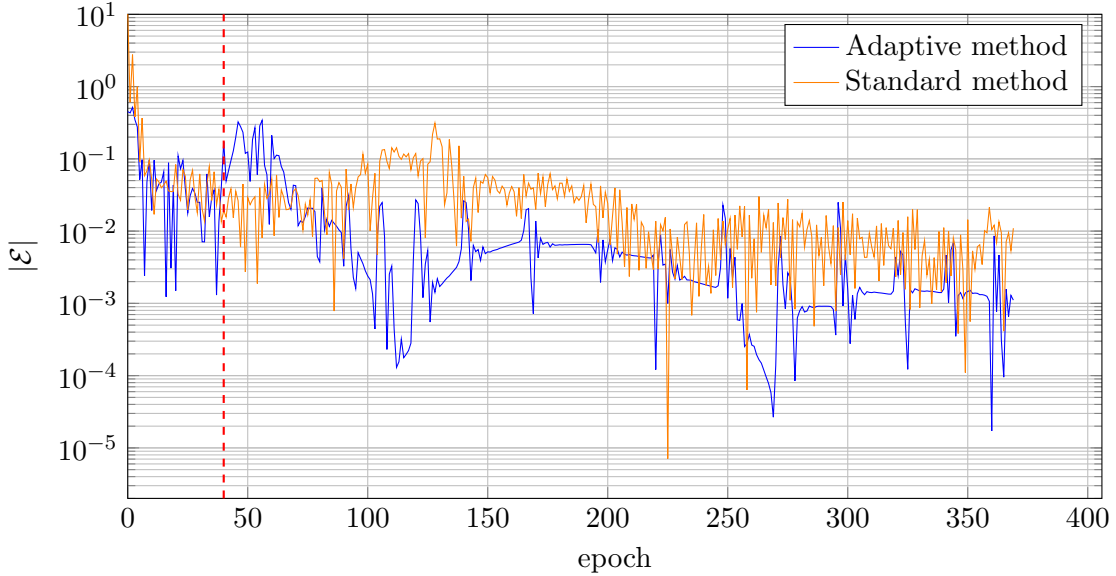
\begin{figure*}

\centering
\begin{tikzpicture}
\begin{axis}[
  width=0.9\linewidth,
  height=8cm,
  xlabel={epoch},
  ylabel={$|{\cal E}|$},
  ymode=log,           % souvent utile pour une loss
  grid=both,
  xmin=0,
  ymin=2.e-6,
  ymax=10
]
\addplot[blue]
table[
  col sep=comma,
  x=epoch,
  y expr={abs(\thisrow{Error})},
]{adaptive_PINN_func_error_plot_ex3.csv};
\addlegendentry{Adaptive method}
\addplot[orange] table[
  col sep=comma,
  x=epoch,
  y expr={abs(\thisrow{Error_standard})},
  mark=none
]{adaptive_PINN_func_error_plot_ex3.csv};
% \node[blue] at (axis cs:900,1e-4) {Standard method};
\addplot[red, thick, dashed]
  coordinates {(40,\pgfkeysvalueof{/pgfplots/ymin})
               (40,\pgfkeysvalueof{/pgfplots/ymax})};
\addlegendentry{Standard method}
\end{axis}
\end{tikzpicture}
\caption{Case 3: Absolute value of the functional error 
$|\mathcal{E}(u)| = |J(u) - J(u_{\mathcal{N}})|$ during the training. The plotted values correspond to the error averaged over every 10 epochs. The vertical dashed red line indicates the iteration at which resampling is performed.}
\label{fig:case3}
\end{figure*}

The results are shown in Figure~\ref{fig:case3}. Here, the networks produce solutions $u_\theta$ with highly oscillatory behavior. Although we show the averages of the functional error over 10 sucessive epochs, the results are still noisy. Nevertheless, the adaptive procedure produces results which a superior by about one order of magnitude.

\subsection{Application to Deep Ritz}  

For applying the estimator based adaptive sampling to Deep Ritz, we have to be aware of the special role of the loss function~\eqref{DR:loss}, which is an approximation to the continuous energy~\eqref{loss:deepritz}. In contrast to the classical PINN loss~\eqref{PINN:loss}, which is just penalizing the residual, Deep Ritz is based on an approximation of the energy functional itself. When choosing non-uniform collocation points, quadrature weights should be added to maintain a proper integral of the domain $\Omega$. 

\subsubsection{Case 4: Simple Poisson  problem}
On the square $\Omega = (0,\pi)^2$ we consider the  problem
\[
                 -\Delta u = 2 \sin(x) \sin(y)\text{ in }\Omega,\quad
                 u=0\text{ on }\partial\Omega.
\]
The error functional is the average over an inscribed disc
\[
C:=\{x\in\mathbb{R}^2\,:\, \left(x_1-\frac{\pi}{2}\right)^2+\left(x_2-\frac{\pi}{2}\right)^2 \le 1\}
\]
giving rise to the adjoint problem
\[
-\Delta z = \mathds{1}_\mathcal{C}\text{ in }\Omega,\quad
z=0\text{ on }\partial\Omega
\]
and the functional $J$ 
\[
  J(u) = \int_\mathcal{C} u(x)\,\text{d}x.
\] 

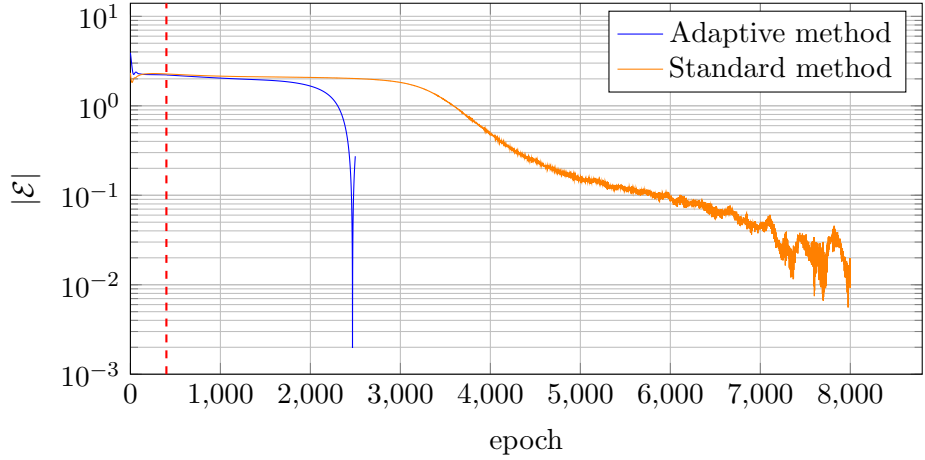
\begin{figure*}[htbp]
\centering
\begin{tikzpicture}
\begin{axis}[
   width=12cm,
  height=6.5cm,
  xlabel={epoch},
  ylabel={$|{\cal E}|$},
  ymode=log,           % souvent utile pour une loss
  grid=both,
  xmin=0,
  ymin=1.e-3,
  ymax=14
]
  \addplot[blue]
table[
  col sep=comma,
  x=epoch,
  y expr={abs(\thisrow{Error})},
]{adaptive_DR_func_plot_ex1.csv};
\addlegendentry{Adaptive method}

\addplot[orange] table[
  col sep=comma,
  x=epoch,
  y expr={abs(\thisrow{Functional error_standard})},
  mark=none
]{adaptive_DR_func_standard_plot_ex1.csv};
 \node[blue] at (axis cs:900,1e-4) {Standard method};
\addplot[red, thick, dashed]
  coordinates {(400,\pgfkeysvalueof{/pgfplots/ymin})
               (400,\pgfkeysvalueof{/pgfplots/ymax})};
\addlegendentry{Standard method}
\end{axis}
\end{tikzpicture}
\caption{Case 4: Absolute value of the functional error 
$|\mathcal{E}(u)| = |J(u) - J(u_{\mathcal{N}})|$ during the training. The vertical dashed red line indicates the iteration at which resampling is performed.}\label{fig:case4}
\end{figure*}

Figure~\ref{fig:case4} shows the slope of the functional error for uniform and adaptive sampling. Here, the effect of error estimator based sampling is even bigger than for the classical PINN formulation. 

\subsubsection{ Case 5: Low dimensional problem with non-smooth functional}

\begin{figure*}[htbp]
\centering
\begin{tikzpicture}
\begin{axis}[
   width=12cm,
  height=6.5cm,
  xlabel={epoch},
  ylabel={$|{\cal E}|$},
  ymode=log,           % souvent utile pour une loss
  grid=both,
  xmin=0,
  ymin=2.1e-3,
  ymax=10
]
\addplot[blue]
table[
  col sep=comma,
  x=epoch,
  y expr={abs(\thisrow{Error})},
]{adaptive_DR_func_plot_ex2.csv};
\addlegendentry{Adaptive method}
\addplot[orange] table[
  col sep=comma,
  x=epoch,
  y expr={abs(\thisrow{Functional error_standard})},
  mark=none
]{adaptive_DR_func_standard_plot_ex2.csv};
% \node[blue] at (axis cs:900,1e-4) {Standard method};
\addplot[red, thick, dashed]
  coordinates {(400,\pgfkeysvalueof{/pgfplots/ymin})
               (400,\pgfkeysvalueof{/pgfplots/ymax})};
\addlegendentry{Standard method}
\end{axis}
\end{tikzpicture}
\caption{Case 5: Absolute value of the functional error 
$|\mathcal{E}(u)| = |J(u) - J(u_{\mathcal{N}})|$ during the training. The vertical dashed red line indicates the iteration at which resampling is performed.}
\label{fig:case5}
\end{figure*}
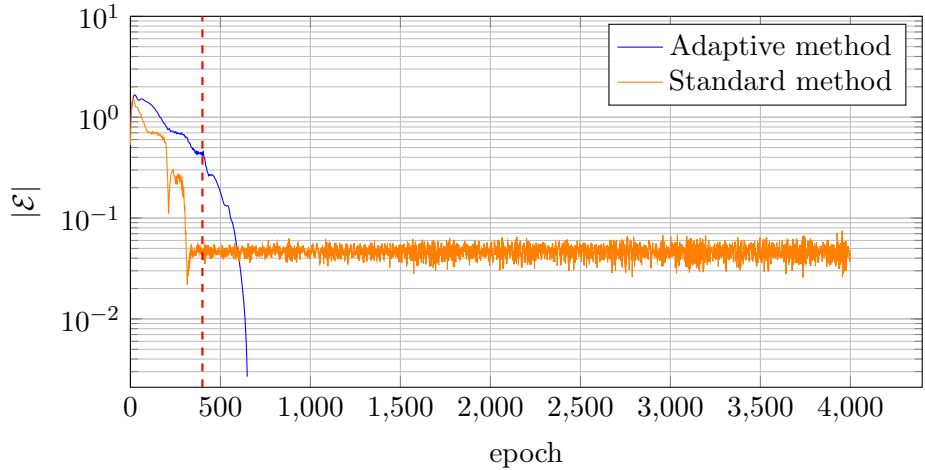

We  consider 
   \begin{equation}\left\{
       \begin{aligned}
                 -\Delta u &= \frac{\pi^2}{2} \cos\left( \frac{\pi}{2} x \right) \cos\left( \frac{\pi}{2} y \right) , \quad (x,y) \in \Omega = ]-1,1[^2 \\
                  u|_{_{\partial \Omega}}& = 0 
       \end{aligned}\right.
   \end{equation} and the adjoint equation 
   \begin{equation}\left\{
       \begin{aligned}
                 -\Delta z &= 1 , \quad (x,y) \in \Omega = [-1,1]^2 \\
                  z|_{_{\partial \Omega}}& = 0,
       \end{aligned}\right.
   \end{equation}
   Then the functional $J$ is defined by  
   \begin{equation*}
       J(u) = \int_\Omega u(x) dx .
   \end{equation*}
   We use a fully residual neural network composed of 2 blocks, each consisting
of two fully connected layers, with 32 neurons per layer, to approximate the exact
solution of the equations. The neural networks are trained for 1500 epochs. For
the adaptive method, the training is first performed using randomly selected
points. At epoch 400, we compute the functional error estimator we resample
the interior points and use them to continue the training of the neural network. Figure \ref{fig:case5} shows the convergence of the functional error. 

\begin{comment}
   \textbf{Example 3 :}
   We consider 
   \begin{equation}\left\{
       \begin{aligned}
                 -\Delta u &=1 , \quad (x,y) \in \Omega = (]0,2[ \times ]0,3[ ) - ( ]0,1[\times ]1,3[)  \\
                  u_{_{\partial \Omega}}& = 0 
       \end{aligned}\right.
   \end{equation} and the adjoint-equation 
   \begin{equation}\left\{
       \begin{aligned}
                 \Delta u &=1 , \ \quad (x,y) \in \Omega = ([0,2] \times [0,3] ) - ( [0,1]\times [1,3])  \\
                  u_{_{\partial \Omega}}& = 0,
       \end{aligned}\right.
   \end{equation}
   Then the functional $J$ is defined by
   \begin{equation*}
       J(u) = \int_{\Omega} u(x) dx .
   \end{equation*}
   \includegraphics[scale=0.8]{Example3.png}
\end{comment} 

\newpage 
\section{Conclusion}

We have derived a fully PINN-based a posteriori error estimator that allows to approximate errors in functional outputs during training and that further can be used for adaptive sampling of collocation points. The most important observation is that the Dual Weighted Residual method can easily be applied to all types of discretizations, if just the underlying problem itself allows for a variational formulation. The lacking linear structure of neural network sets means that concepts like Galerkin orthogonality do not hold. The DWR estimator benefits from this, as no orthogonality constraints must be considered for approximating the weights. All numerical examples demonstrate very good estimator effectivities from early in the training process, although we cannot expect that primal and adjoint solution are already well approximated. 

To localize the error we formulated the discrete problems as variations over the network parameters $\theta$ instead of the network functions $u_\theta$. This gives rise to the definition of the derivative sets $V_{\cal N}'(u_\theta)$ to be used to localize the weights. Hereby, we could improve the quality of the refinement indices and use them for sampling collocation points. Several numerical examples show that adaptive sampling significantly accelerates the training and also helps to reach smaller functional errors of the trained networks. 

Several open questions remain for future research. From the analytical perspective the most intriguing question is the role of orthogonality and the reason why very rough approximations of primal and adjoint problems already yield highly accurate estimators, see also~\cite{RothSchroederWick2022,MR}.

Technically, it remains to improve the actual sampling of the collocation points. Although adaptive sampling is successful in all test cases, we believe that the theoretical potential should be even higher. Furthermore, in the context of the Deep Ritz method, the objective is to reformulate the loss functional using the estimated measure, so that the weights of the sampled points explicitly appear in the loss. This procedure may lead to the construction of an approximation of the estimated measure itself. To this end, we aim to draw inspiration from a recent work by \cite{Wan25}, in which the Knothe–Rosenblatt rearrangement is used to construct such an approximation.

\section*{Acknowledgements}
The authors acknowledge support of the \emph{Deutsche Forschungsgemeinschaft (DFG)} by means of the Research Training Group \emph{MathCore} (grant number 314838170). The authors further acknowledge support of the \emph{Deutsche Forschungsgemeinschaft (DFG)} under the grant with number 537063406. Furthermore, the authors acknowledge support of the \emph{German Federal Ministry of Research Technology and Space} under the grant with number 03VP13242.

% \appendix
% \section{My Appendix}
% Appendix sections are coded under \verb+\appendix+.

% \verb+\printcredits+ command is used after appendix sections to list 
% author credit taxonomy contribution roles tagged using \verb+\credit+ 
% in frontmatter.


\bibliographystyle{plain}

\bibliography{ml_references}

\end{document}